\documentclass[12pt,oneside]{article}

\usepackage{amsmath}
\usepackage{amssymb}
\usepackage{color}
\usepackage{rotating}

\usepackage{graphicx}
\usepackage{subfigure}

\usepackage{url}

\usepackage{appendix}

\author{Hui Zhou\footnote{E-mail addresses:\newline zhouhpku17@pku.edu.cn, huizhou@math.pku.edu.cn, zhouhlzu06@126.com.}\\ \footnotesize{School of Mathematical Sciences, Peking University, Beijing, 100871, P.~R.~China.} }

\title{\Large A graph theoretic characterization of the classical generalized hexagon on $364$ vertices}

\def\Aut{{\sf Aut}}

\def\PSL{{\sf PSL}}
\def\PSU{{\sf PSU}}
\def\G{{\sf G}}
\def\R{{\sf R}}

\def\diam{{\rm diam}}
\def\deg{{\rm deg}}

\def\A{{\sf A}}
\def\S{{\sf S}}
\def\GL{{\sf GL}}
\def\Z{{\mathbb{Z}}}

\newtheorem{theorem}{Theorem}

\newtheorem{lemma}[theorem]{Lemma}

\newcommand*{\QEDA}{\hfill\ensuremath{\blacksquare}}  

\begin{document}

\maketitle

\begin{abstract}


A tetravalent $2$-arc-transitive graph of order $728$ is either the known $7$-arc-transitive incidence graph of the classical generalized hexagon $GH(3,3)$ or a normal cover of a $2$-transitive graph of order $182$ denoted $A[182,1]$ or $A[182,2]$ in the $2009$ list of Poto\v{c}nik.
\end{abstract}

%


%

\section{Preliminaries}\label{sec priminaries}

In this note, all graphs are finite, simple, connected and undirected. An ordered pair of adjacent vertices is called an arc. Let $\Gamma$ be a graph. We use $V\Gamma$, $E\Gamma$, $A\Gamma$ and $\Aut(\Gamma)$ to denote the vertex-set, the edge-set, the arc-set and the full automorphism group of $\Gamma$, respectively. The distance between two vertices $u$ and $v$ of $\Gamma$, denoted by $\partial_\Gamma(u,v)$, is the length of a shortest path connecting $u$ and $v$ in $\Gamma$. The diameter of $\Gamma$, denoted by $\diam(\Gamma)$, is the maximum distance occurring over all pairs of vertices. Fix a vertex $v\in V\Gamma$. For $0\leqslant i\leqslant \diam(\Gamma)$, we use $\Gamma_i(u)$ to denote the set of vertex $u$ with $\partial_\Gamma(u,v)=i$. For convenience, we usually use $\Gamma(v)$ to denote $\Gamma_1(v)$. The degree of $v$, denoted by $\deg_\Gamma(v)$ or simply $\deg(v)$, is the number of vertices adjacent to $v$ in $\Gamma$, i.e. $\deg_\Gamma(v)=|\Gamma(v)|$. The graph $\Gamma$ is called regular with valency $k$ (or $k$-regular) if the degree of each vertex of $\Gamma$ is $k$. The girth of $\Gamma$ is the length of a shortest cycle of $\Gamma$.

Let $G\leqslant \Aut(\Gamma)$. The graph $\Gamma$ is called $G$-vertex-transitive ($G$-arc-transitive, respectively), if $G$ is transitive on the vertex-set $V\Gamma$ (the arc-set $A\Gamma$, respectively). Let $s\geqslant 1$. An $s$-arc of $\Gamma$ is an $(s+1)$-tuple of vertices of $\Gamma$ in which every two consecutive vertices are adjacent and every three consecutive vertices are pairwise distinct. The graph $\Gamma$ is called $(G,s)$-arc-transitive if it is $G$-vertex-transitive and $G$ is also transitive on $s$-arcs of $\Gamma$. The graph $\Gamma$ is called $(G,s)$-transitive if it is $(G,s)$-arc-transitive but not $(G,s+1)$-arc-transitive. The graph $\Gamma$ is called $(G,s)$-distance-transitive, if for each $1\leqslant i\leqslant s$ the group $G$ is transitive on the orderer pairs of form $(u,v)$ with $\partial_\Gamma(u,v)=i$. A $(G,s)$-distance-transitive graph is called $G$-distance-transitive (or distance-transitive), if $s$ is the diameter of the graph (and $G$ is the full automorphism group of the graph). By definitions, $s$-arc-transitivity implies $s$-distance-transitivity.

Take a vertex $\alpha$ of a vertex-transitive graph $\Gamma$ and let $0\leqslant i\leqslant \diam(\Gamma)$ , then the size $\kappa_i=|\Gamma_i(\alpha)|$ is independent of $\alpha$.

Let $\Gamma$ be a $(G,s)$-distance-transitive graph with valency $k\geqslant 3$. For $1\leqslant i\leqslant s$, take $\alpha\in V\Gamma$ and $\beta\in \Gamma_i(\alpha)$, then the intersection numbers
\begin{eqnarray*}
c_i&=&|\Gamma_{i-1}(\alpha)\cap \Gamma(\beta)|,\\
a_i&=&|\Gamma_i(\alpha)\cap \Gamma(\beta)|,\\
b_i&=&|\Gamma_{i+1}(\alpha)\cap \Gamma(\beta)|
\end{eqnarray*}
are independent of $\alpha$ and $\beta$, and $c_i+a_i+b_i=k$. The $s$-partial intersection array of $\Gamma$ is
\begin{equation*}
\iota(\Gamma,s)=\left\{\begin{array}{ccccc}
               * & c_1 & c_2 & \cdots & c_s \\
               0 & a_1 & a_2 & \cdots & a_s \\
               k & b_1 & b_2 & \cdots & b_s
               \end{array}\right\}.
\end{equation*}
Then $c_1=1$. For convenience, we let $b_0=k$. Take $\alpha\in V\Gamma$ and let $1\leqslant i\leqslant s$. By counting edges between $\Gamma_{i-1}(\alpha)$ and $\Gamma_i(\alpha)$, we have \begin{center}$\kappa_{i-1} b_{i-1}=\kappa_i c_i$.\end{center}

From~\cite[Table~2.4 on pages 135-136]{GorensteinCFSG} or~\cite[Chapter~7]{CameronPermutationGroup}, one may obtain the following lemma by checking the orders of non-abelian finite simple groups.

\begin{lemma}\label{lem non-solvable minimal normal subgroup N}
Let $N$ be a non-abelian finite simple group with $|\pi(N)|\geqslant 3$ and $\pi(N)\subseteq \{2,3,7,13\}$. Then $N$ is one of the groups listed in Table~\ref{table non-abelian finite simple 2-3-7-13-group}.
\end{lemma}

\begin{table}[!ht]
\centering
    \begin{tabular}{|l|rcl|}
        \hline
        group & & order & \\
        \hline
        $\PSL(3,2)=\PSL(2,7)$ & $168$ & $=$ & $2^3\cdot 3\cdot 7$\\
        $\PSL(2,8)$ & $504$ & $=$ & $2^3\cdot 3^2\cdot 7$\\
        $\R(3)={}^{2}\G_2(3)=\PSL(2,8):3$ & $1512$ & $=$ & $2^3\cdot 3^3\cdot 7$\\
        $\PSL(3,3)$ & $5616$ & $=$ & $2^4\cdot 3^3\cdot 13$\\
        $\PSU(3,3)={}^{2}\A_2(3)$ & $6048$ & $=$ & $2^5\cdot 3^3\cdot 7$\\
        $\G_2(2)=\PSU(3,3):2$ & $12096$ & $=$ & $2^5\cdot 3^3\cdot 7$\\
        \hline
        $\PSL(2,27)$ & $9828$ & $=$ & $2^2\cdot 3^3\cdot 7\cdot 13$\\
        $\G_2(3)$ & $4245696$ & $=$ & $2^6\cdot 3^6\cdot 7\cdot 13$\\
        \hline
    \end{tabular}
    \caption{Non-abelian finite simple $\{2,3,7,13\}$-group.}\label{table non-abelian finite simple 2-3-7-13-group}
\end{table}


\section{The Proof}\label{section Proof}

Let $\Gamma$ be a tetravalent $(G,2)$-arc-transitive graph of order \begin{center}$|V\Gamma|=728=8pq$\end{center} where $G=\Aut(\Gamma)$ and $(p,q)=(7,13)$. Then $\Gamma$ is $(G,s)$-transitive for some $s\geqslant 2$. Take a vertex $\alpha\in V\Gamma$. Then the vertex stabilizer $G_\alpha$ and $s$ are listed in Table~\ref{table stabilizer of tetravalent 2-arc-transitive graphs}~\cite[Lemma~2.6]{LiLuWangTetraSquarefreeOrder}. Hence $|G_\alpha|~\bigl|~ 2^4\cdot 3^6$. Since $\Gamma$ is $(G,2)$-arc-transitive, we have $G_\alpha^{\Gamma(\alpha)}\lesssim \S_4$ is $2$-transitive on $\Gamma(\alpha)$, i.e. $G_\alpha^{\Gamma(\alpha)}\cong \A_4$ or $\S_4$. Then $12~\bigl|~|G_\alpha^{\Gamma(\alpha)}|~\bigl|~|G_\alpha|$, and so $|G_\alpha|=2^{i'}\cdot 3^{j'}$ where $2\leqslant i'\leqslant 4$ and $1\leqslant j'\leqslant 6$. Since $\Gamma$ is $G$-vertex-transitive, by Frattini argument on permutation groups, we have \begin{center}$|G|=|G_\alpha|\cdot |V\Gamma|=2^i\cdot 3^j\cdot p\cdot q=2^i\cdot 3^j\cdot 7\cdot 13\leqslant 8491392$\end{center} where $5\leqslant i=3+i'\leqslant 7$ and $1\leqslant j=j'\leqslant 6$. Note that $2$-transitivity implies primitivity. We have that the permutation group $G_\alpha^{\Gamma(\alpha)}$ is primitive on $\Gamma(\alpha)$.

\begin{table}[!ht]
\centering
    \begin{tabular}{|c|c|c|c|c|}
        \hline
        $s$ & $2$ & $3$ & $4$ & $7$ \\
        \hline
        $G_\alpha$ & \begin{tabular}{c}$\A_4$ or\\ $\S_4$\end{tabular} & \begin{tabular}{c}$\Z_3\times \A_4$,\\ $(\Z_3\times \A_4).\Z_2$ or\\ $\S_3\times \S_4$\end{tabular} & $\Z_3^2:\GL(2,3)$ & $[3^5]:\GL(2,3)$ \\
        \hline
        $|G_\alpha|$ & \begin{tabular}{c}$12=2^2\cdot 3$ or\\ $24=2^3\cdot 3$\end{tabular} & \begin{tabular}{c}$36=2^2\cdot 3^2$,\\ $72=2^3\cdot 3^2$ or\\ $144=2^4\cdot 3^2$\end{tabular} & $432=2^4\cdot 3^3$ & \begin{tabular}{c}$11664$\\ $=2^4\cdot 3^6$\end{tabular} \\
        \hline
    \end{tabular}
    \caption{The stabilizer $G_\alpha$ of a tetravalent $(G,s)$-transitive graph}\label{table stabilizer of tetravalent 2-arc-transitive graphs}
\end{table}

First, we suppose $G$ has a solvable minimal normal subgroup $N$. Then $N$ is an elementary $r$-group with $r\in \{2,3,p,q\}$. Suppose $|N|=r^e$ for some $e\geqslant 1$. Note that $N\neq N_\alpha$. Let $|N_\alpha|=r^{e'}$. Then $0\leqslant e'<e$. By Frattini argument, we have $|N|=|N_\alpha|\cdot |\alpha^N|$. Hence $|\alpha^N|=r^{e-e'}$. Consider the quotient graph $\Gamma_N$. We have $|V\Gamma|=|\alpha^N|\cdot |V\Gamma_N|$, which means $r^{e-e'}~\bigl|~|V\Gamma|=8pq$. So $r\neq 3$, i.e. $r\in \{2,p,q\}$. The order of the quotient graph $\Gamma_N$ is $|V\Gamma_N|=\frac{|V\Gamma|}{|\alpha^N|}=\frac{8pq}{r^{e-e'}}>2$. So $N$ is semi-regular on $V\Gamma$, $V\Gamma_N$ is a tetravalent $(G/N,2)$-arc-transitive graph and $\Gamma$ is a cover of $\Gamma_N$~\cite[Theorem~4.1]{Preager-O'Nan-Scott-2AT}. Now $N=\Z_p$, $\Z_q$, $\Z_2$, $\Z_2^2$ or $\Z_2^3$. Then $\Gamma_N$ is a $2$-arc-transitive graph of order in $\{8p,8q,4pq,2pq,pq\}=\{56,91,104,182,364\}$. By~\cite[Table~3 on pages 1334-1335]{4valent-2at}, we have that $\Gamma_N$ is $A[182,1]$ or $A[182,2]$ from~\cite{4valent-2at}, and so $N=\Z_2^2$. These two graphs  are $2$-transitive, i.e. they are $2$-arc-transitive, but not $s$-arc-transitive for $s\geqslant 3$

%

Now we suppose $G$ has no solvable minimal normal subgroups. Take a minimal normal subgroup $N$ of $G$. Then $N$ is non-solvable. By Burnside's $p^aq^b$-theorem in~\cite[Theorem~4.130 on page 239]{GorensteinCFSG} or~\cite[page~36]{CameronPermutationGroup}, we have $|\pi(N)|\geqslant 3$. Let $N=T^k$ where $T$ is a non-abelian simple group and $k\geqslant 1$. Note that $\pi(T)=\pi(N)\subseteq \pi(G)=\{2,3,p,q\}$. Then $p$ or $q$ is in $\pi(T)$ which implies $k=1$, i.e. $N$ is a non-abelian simple group. So $N$ is one of the groups in Table~\ref{table non-abelian finite simple 2-3-7-13-group}. Since $|N_\alpha|~\bigl|~|G_\alpha|~\bigl|~2^4\cdot 3^6$, for any prime $r\in \pi(N)\setminus \{2,3\}$, we have $r~\bigl|~|\alpha^N|=\frac{|N|}{|N_\alpha|}$. Consider the quotient graph $\Gamma_N$. Note that $|\alpha^N|~\bigl|~\frac{|V\Gamma|}{|V\Gamma_N|}~\bigl|~8pq$. Then $|\alpha^N|_2~\bigl|~8$ and $3\nmid |\alpha^N|$, which implies $|N|_3=|N_\alpha|_3$. By Table~\ref{table non-abelian finite simple 2-3-7-13-group}, we have $|N|_3\neq 1$. This implies $N_\alpha\neq 1$. So $N_\alpha^{\Gamma(\alpha)}\neq 1$ and $N$ is not semi-regular on $V\Gamma$. Note that $N_\alpha^{\Gamma(\alpha)}\trianglelefteq G_\alpha^{\Gamma(\alpha)}$. Thus $N_\alpha^{\Gamma(\alpha)}$ is transitive on $\Gamma(\alpha)$. We have \begin{center}$4~\bigl|~|N_\alpha^{\Gamma(\alpha)}|~\bigl|~|N_\alpha|$.\end{center} If $|V\Gamma_N|\geqslant 3$, then $N$ is semi-regular on $V\Gamma$~\cite[Theorem~4.1]{Preager-O'Nan-Scott-2AT}. So we have $|V\Gamma_N|\leqslant 2$, i.e. $|\alpha^N|=4pq$ or $8pq$. This means \begin{center}$|N|=4pq\cdot |N_\alpha|$ or $8pq\cdot |N_\alpha|$.\end{center}
If $\pi(N)=\{2,3,r\}$ where $r\in\{p,q\}$, then $|\alpha^N|\leqslant 8r\leqslant 8\max\{p,q\}$, and so $|V\Gamma_N|\geqslant \min\{p,q\}\geqslant 3$. This is a contradiction. Hence \begin{center}$\pi(N)=\{2,3,p,q\}$, and so $N=\PSL(2,27)$ or $\G_2(3)$.\end{center}
If $N=\PSL(2,27)$, then $8\nmid |N|$, and so we have $|N|=4pq\cdot |N_\alpha|$ and $|N_\alpha|=3^3$, which is a contradiction. Hence we have $N=\G_2(3)$.


%
%


Let $N=\G_2(3)$. Then $|N|=4pq\cdot |N_\alpha|$ or $8pq\cdot |N_\alpha|$. So $|N_\alpha|=11664=2^4\cdot 3^6$ or $|N_\alpha|=5832=2^3\cdot 3^6$. By Table~\ref{table stabilizer of tetravalent 2-arc-transitive graphs}, we have \begin{center}$(s,G_\alpha)=(7,[3^5]:\GL(2,3))$.\end{center}
Let $|N|=8pq\cdot |N_\alpha|$. Then $|N_\alpha|=2^3\cdot 3^6$, $|\alpha^N|=8pq$ and $|V\Gamma_N|=1$. This implies that $N$ is transitive on $V\Gamma$. By~\cite[Theorem~1.1]{LiSeressSong-sATandNormalSubgroups}, we get that $\Gamma$ is $(N,7)$-transitive. By Table~\ref{table stabilizer of tetravalent 2-arc-transitive graphs}, we have $|N_\alpha|=2^4\cdot 3^6$. This is a contradiction.
Hence \begin{center}$|N|=4pq\cdot |N_\alpha|$, $|\alpha^N|=4pq$ and $|V\Gamma_N|=2$.\end{center} So $\Gamma$ is bipartite. Now $|N_\alpha|=2^4\cdot 3^6=|G_\alpha|$. Since $N_\alpha\leqslant G_\alpha$, we get \begin{center}$N_\alpha=G_\alpha$.\end{center} Note that $|G|=|G_\alpha|\cdot |V\Gamma|=|N_\alpha|\cdot 8pq=2|N|$. Thus we have \begin{center}$G=N.2=\G_2(3).2$.\end{center} Let $g$ be the girth of $\Gamma$. By~\cite[Proposition・17.2 on page 131]{Biggs-AGT}, we get that $g\geqslant 2s-2=12$. Since $\Gamma$ is bipartite, there is no odd cycles in $\Gamma$. So $g$ is even. Let $g=2t+2$. Then $\diam(\Gamma)\geqslant \frac{g}{2}=t+1$, $t\geqslant 5$ and $\Gamma$ is $(G,t)$-arc-transitive. The $t$-partial intersection array of $\Gamma$ is
\begin{equation*}
\iota(\Gamma,t)=\left\{\begin{array}{cccccccc}
               * & 1 & 1 & 1 & 1 & 1 & \cdots & c_t=1 \\
               0 & 0 & 0 & 0 & 0 & 0 & \cdots & a_t=0 \\
               4 & 3 & 3 & 3 & 3 & 3 & \cdots & b_t=3
               \end{array}\right\}.
\end{equation*}
Hence $3^6-1=728=|V\Gamma|\geqslant \sum\limits_{i=0}^t|\Gamma_i(\alpha)|=1+4+4\cdot 3+4\cdot 3^2 + \cdots +4\cdot 3^{t-1}=2\cdot 3^t-1$. This implies $t<6$. So we have $t=5$, $g=12$ and $\diam(\Gamma)\geqslant 6$. Now we suppose $\diam(\Gamma)\geqslant 7$. Then $c_6\leqslant 3=b_5$. We have $|\Gamma_6(\alpha)|=\frac{b_5}{c_6}|\Gamma_5(\alpha)|\geqslant |\Gamma_5(\alpha)|=4\cdot 3^{t-1}$. So $|V\Gamma|\geqslant \sum\limits_{i=0}^6|\Gamma_i(\alpha)|=2\cdot 3^t-1+4\cdot 3^{t-1}=3^6+3^4-1>3^6-1=|V\Gamma|$. This is a contradiction. Hence $\diam(\Gamma)=6$. Then $\Gamma$ is $G$-distance-transitive. By~\cite[page~222]{BCN} or~\cite{4dr}, $\Gamma$ is the incidence graph of the known generalized hexagon $GH(3,3)$, i.e. the generalized dodecagon of order $(1,3)$. We denote this graph by $\Gamma_7$. The graph $\Gamma_7$ is actually $7$-arc-transitive by~\cite[last prargraph on the first page]{Gardiner-ATprimevalency} or~\cite[17g on page 137]{Biggs-AGT}, and it is $7$-transitive with smallest order.

The construction of $\Gamma_7$ can be found in~\cite{Bensen1966}. Another construction of $\Gamma_7$ is by orbital graph method. The full automorphism group of $\Gamma_7$ is $\G_2(3).2$ with vertex stabilizer of order $11664$. The almost simple group $\G_2(3).2$ with socle $\G_2(3)$ is in the database of almost simple groups of Magma. The group $\G_2(3).2$ has only one conjugate class of subgroups of order $11664$. The graph $\Gamma_7$ is isomorphic to the orbital graph with arc-transitive group $\G_2(3).2$ and an orbit of length four under the action of the stabilizer, and so it can be constructed by Magma. The following is the corresponding Magma code of the construction of $\Gamma_7$.

\begin{verbatim}
D:=AlmostSimpleGroupDatabase();
NumberOfGroups(D,4245696,8491392);
G:=GroupData(D,203)`permrep;
M:=Subgroups(G:OrderEqual:=11664);
#M;
H:=M[1]`subgroup;
GonH:=CosetImage(G,H);
newH:=Stabilizer(GonH,1);
Orb:=Orbits(newH);
for o in Orb do
#o;
end for;
Neb:=Set(Orb[2]);
gamma:=OrbitalGraph(GonH,1,Neb);
Aut:=AutomorphismGroup(gamma);
f,_:=IsIsomorphic(G,Aut);
f;
gamma;
\end{verbatim}

The result graph $\Gamma_7$ with vertices numbered $\{1,2,\ldots,728\}$ and neighbours for each vertex is in Appendix~\ref{appendix the graph gamma7}. Note that, the labels of the vertices are randomly changed.




\section{Acknowledgement}

The author would like to thank Binzhou Xia for his critical reading and valuable comments and suggestions which improve the original manuscript.


%
%



\begin{appendices}

\section{The graph $\Gamma_7$}\label{appendix the graph gamma7}

\begin{figure}[htb]
\centering
\subfigure{\includegraphics[width=0.24\textwidth]{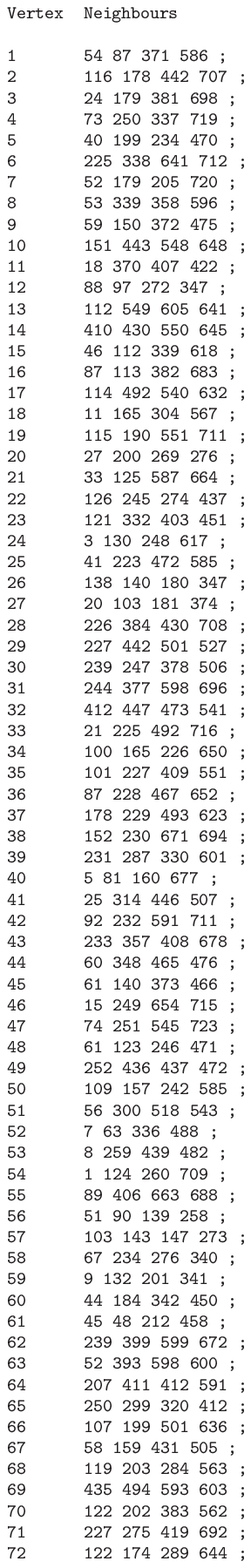}}
\subfigure{\includegraphics[width=0.24\textwidth]{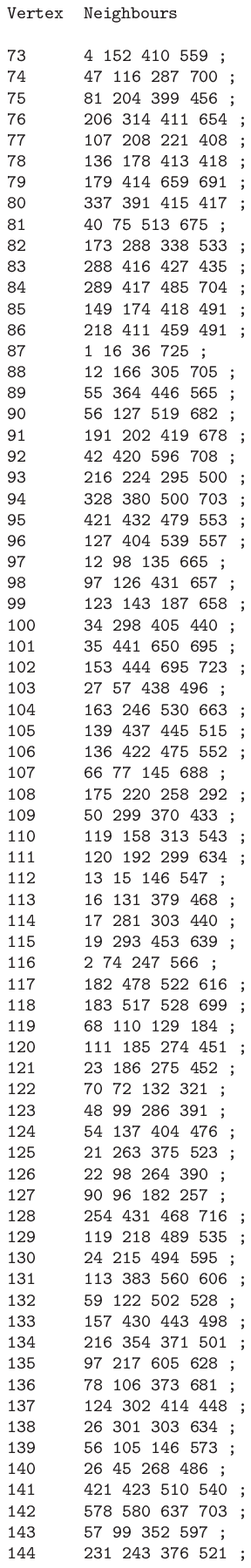}}
\subfigure{\includegraphics[width=0.24\textwidth]{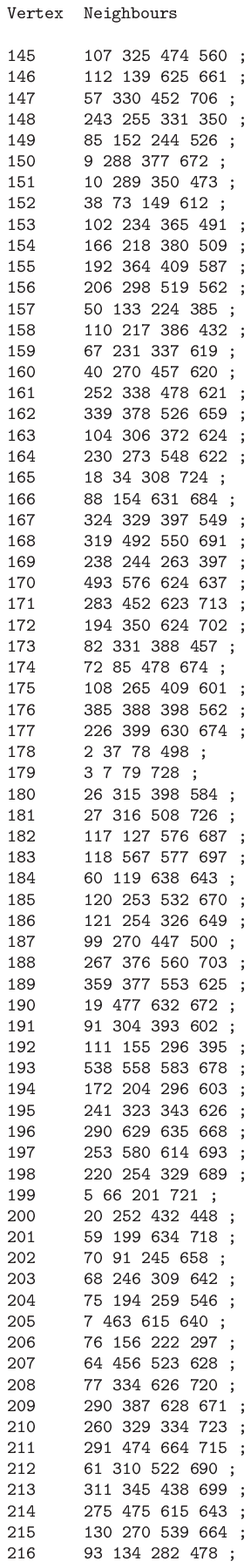}}
\end{figure}

\begin{figure}[htb]
\centering
\subfigure{\includegraphics[width=0.24\textwidth]{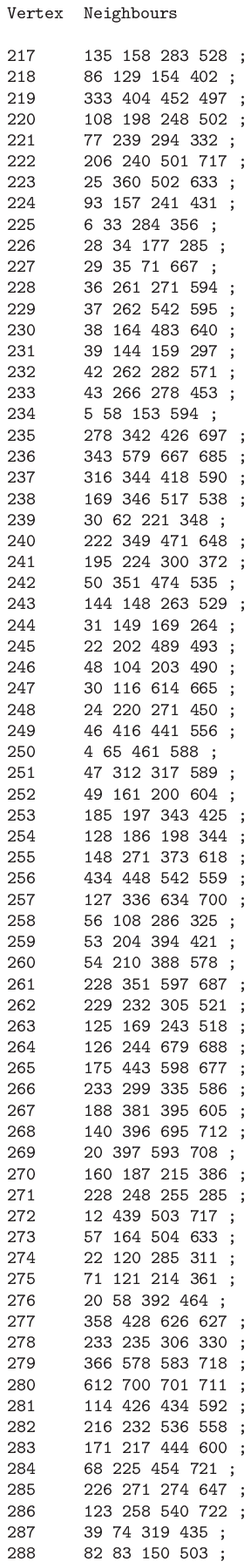}}
\subfigure{\includegraphics[width=0.24\textwidth]{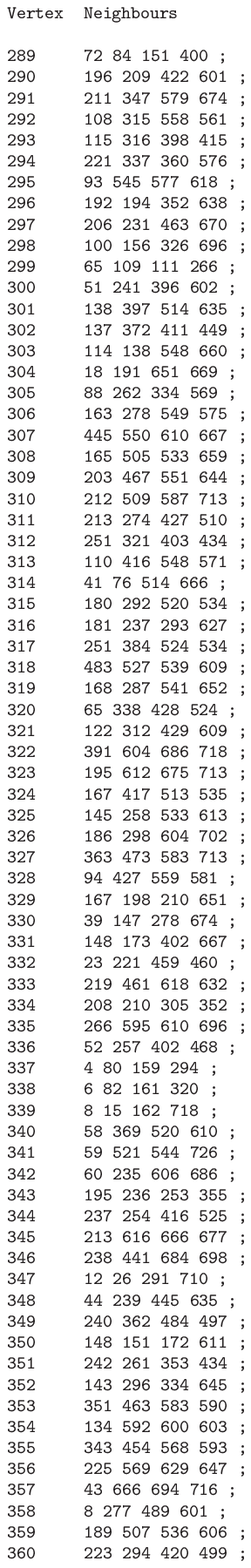}}
\subfigure{\includegraphics[width=0.24\textwidth]{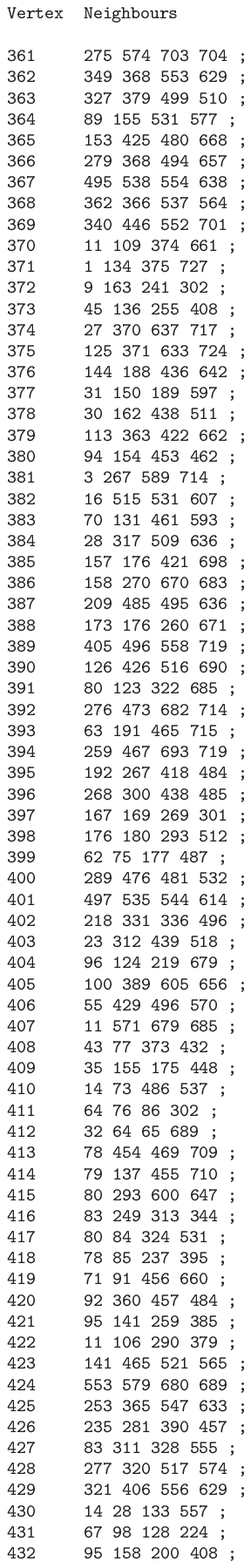}}
\end{figure}

\begin{figure}[htb]
\centering
\subfigure{\includegraphics[width=0.24\textwidth]{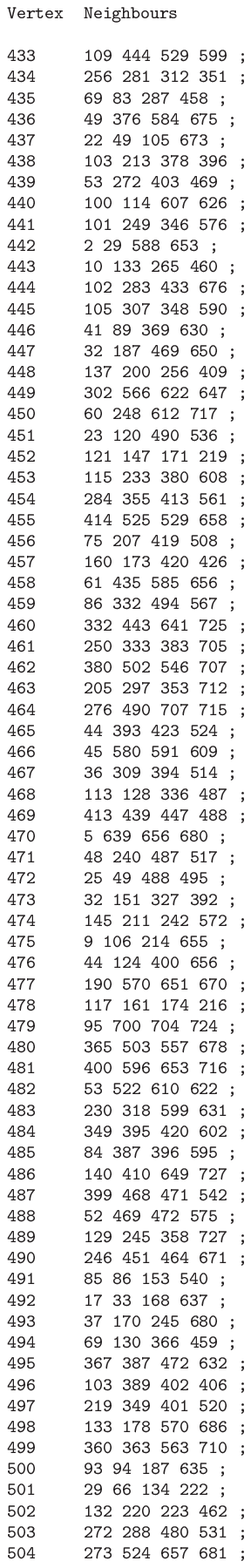}}
\subfigure{\includegraphics[width=0.24\textwidth]{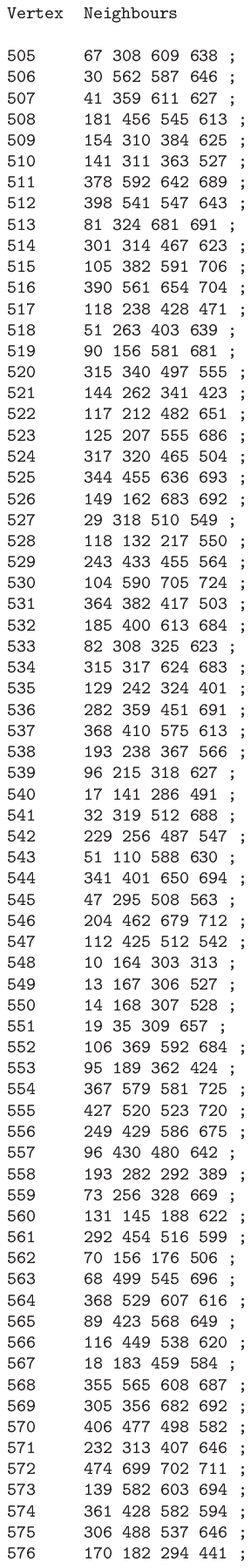}}
\subfigure{\includegraphics[width=0.24\textwidth]{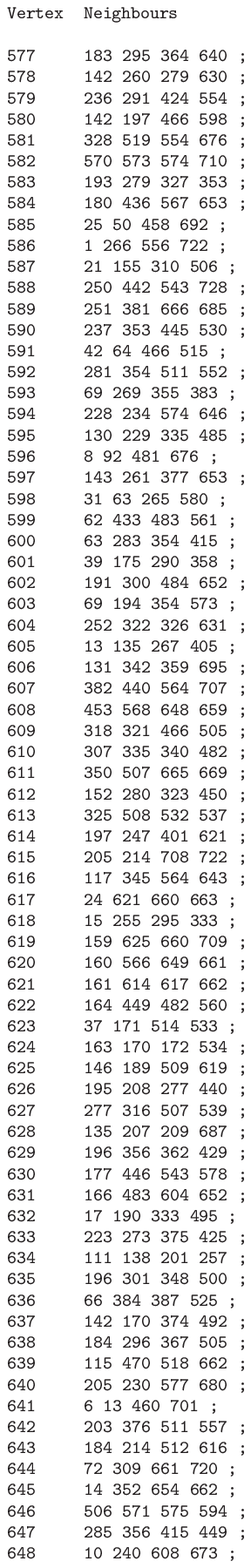}}
\end{figure}

\begin{figure}[htb]
\centering
\subfigure{\includegraphics[width=0.24\textwidth]{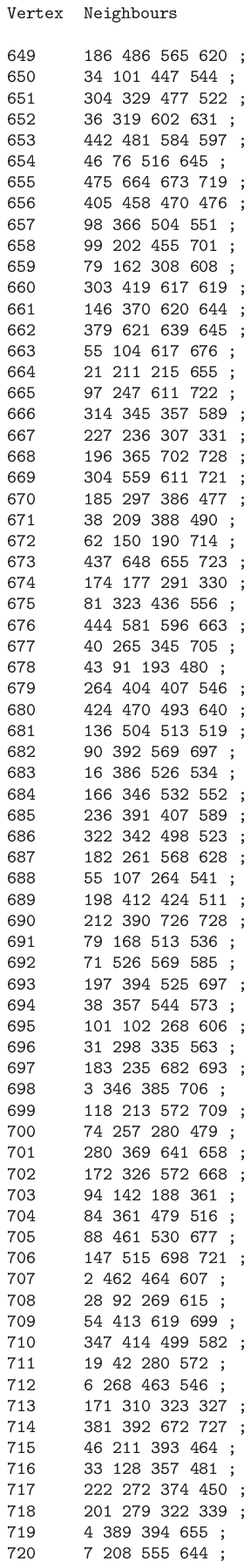}}
\subfigure{\includegraphics[width=0.24\textwidth]{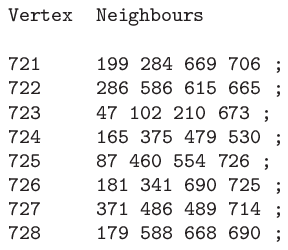}}
\end{figure}

\end{appendices}

\end{document}